\documentclass[aps, prd, floats, floatfix, onecolumn, superscriptaddress, nofootinbib, letterpaper, notitlepage, preprintnumbers]{revtex4-1}

\usepackage{graphicx}
\usepackage{amsmath}
\usepackage{amsfonts}
\usepackage{amssymb}
\usepackage{amstext}
\usepackage{tensor}
\usepackage{fullpage}
\usepackage[english]{babel}
\usepackage{blindtext}
\usepackage{helvet}
\usepackage{microtype}
\usepackage[pdftex,
  pdftitle={Indefinite Integrals of Spherical Bessel Functions},
  pdfauthor={Jolyon K. Bloomfield, Stephen H. P. Face, Zander Moss},
  bookmarks,bookmarksopen=false,
  pdfstartview={FitH},
  linktocpage=true
]{hyperref}

\newcommand{\ci}{\mathrm{Ci}}
\newcommand{\si}{\mathrm{Si}}

\begin{document}

\title{Indefinite Integrals of Spherical Bessel Functions}

\author{Jolyon K. Bloomfield}
\email{jolyon@mit.edu}
\author{Stephen H. P. Face}
\email{face@mit.edu}
\author{Zander Moss}
\email{zander@mit.edu}
\affiliation{Center for Theoretical Physics, Laboratory for Nuclear Science, and Department of Physics, Massachusetts Institute of Technology, Cambridge, MA 02139, USA}

\date{\today}

\preprint{MIT-CTP/4872}

\begin{abstract}
	Highly oscillatory integrals, such as those involving Bessel functions, are best evaluated analytically as much as possible, as numerical errors can be difficult to control. We investigate indefinite integrals involving monomials in $x$ multiplying one or two spherical Bessel functions of the first kind $j_l(x)$ with integer order $l$. Closed-form solutions are presented where possible, and recursion relations are developed that are guaranteed to reduce all integrals in this class to closed-form solutions. These results allow for definite integrals over spherical Bessel functions to be computed quickly and accurately. For completeness, we also present our results in terms of ordinary Bessel functions, but in general, the recursion relations do not terminate.
\end{abstract}

\maketitle

\section{Introduction}

In a model of early universe physics, we encountered integrals of the form
\begin{align}
\int dx \, f(x) j_l(\alpha x) \quad \text{and} \quad \int dx \, f(x) j_k(\alpha x) j_l(\beta x)
\end{align}
when attempting to compute the covariance between the spherical harmonic modes of a field. Here, $j_l(x)$ is the spherical Bessel function of the first kind, and we were interested in integers $k, l \ge 0$, and reals $\alpha, \beta$. Our functions $f(x)$ were slowly-varying functions of $x$ that had been numerically computed at discrete values. As a piecewise-polynomial interpolation of $f(x)$ is an appropriate description of such a function, it was natural to consider indefinite integrals of the form
\begin{align} \label{eq:interest}
\int dx \, x^n j_l(\alpha x) \quad \text{and} \quad \int dx \, x^n j_k(\alpha x) j_l(\beta x)
\end{align}
to numerically compute our integrals. Because the spherical Bessel functions can oscillate wildly when $\alpha x$ and/or $\beta x$ are large, it was important that these integrals be computed analytically.

We searched the literature for integration methods for integrands involving spherical Bessel functions. We found a number of results regarding definite integrals over an infinite range \cite{Mehrem1991, Maximon1991, Mehrem2009, Fabrikant2013}, but very little regarding indefinite integrals. Results for Bessel functions (with integer order) appear to be much more commonplace \cite{Rosenheinrich2016}. NIST's Digital Library of Mathematical Functions \cite{NIST:DLMF} was unhelpful, and even the venerable Gradshteyn and Rhyzik \cite{Gradshteyn} provided little assistance.

Wolfram Research \cite{wolfram-besselj} had some of the integrals that we wanted and Mathematica could integrate a few more, but essentially all these results were written in terms of generalized hypergeometric $\tensor[_p]{F}{_q}$ functions. As these general functions are not available in typical scientific libraries, we set out to derive our own results.

This paper presents analytic results for indefinite integrals of the form \eqref{eq:interest} with integer $n$, integers $k, l \ge 0$ and reals $\alpha$ and $\beta$\footnote{A number of our results likely apply more generally, but we do not investigate this.}. The results have been constructed with an eye towards numerical implementation in a language like C or Python, and hence are written only in terms of special functions that are commonly available in scientific libraries. For some expressions, we find closed-form results, while for others, we derive recursion relations that decrease the the order of spherical Bessel functions appearing in the integrals. For the class of integrals we consider, all recursion relations terminate with closed-form solutions. This only occurs for spherical Bessel functions with integer order, as integrals involving $j_0(x) = \mathrm{sinc}(x)$ are analytically tractable (unlike integrals involving $J_0(x)$). To our knowledge, the recursion relations we derive are new, as are some of the closed-form solutions that we present.

We begin this paper by describing well-known properties of spherical Bessel functions. We then construct some preliminary quantities which appear in the general integrals, before treating the single and product integrals in turn. An appendix translates applicable results to ordinary Bessel functions.

\section{Overview of Spherical Bessel Functions}

We begin by presenting well-known results regarding spherical Bessel functions. Spherical Bessel functions are the solutions to the differential equation
\begin{align}
x^2 \frac{d^2 y}{dx^2} + 2x \frac{dy}{dx} + [x^2 - l(l+1)] y = 0.
\end{align}
They are indexed by the order $l$. Two independent solutions are $j_l(x)$ and $y_l(x)$, the spherical Bessel functions of the first and second kind, respectively. These functions are related to ordinary Bessel functions by
\begin{align}
j_l(x) = \sqrt{\frac{\pi}{2x}} J_{l + 1/2}(x)
\quad \text{and} \quad
y_l(x) = \sqrt{\frac{\pi}{2x}} Y_{l + 1/2}(x) = (-1)^{l+1} \sqrt{\frac{\pi}{2x}} J_{- l - 1/2}(x).
\end{align}
The spherical Bessels can be computed from Rayleigh's formulas
\begin{align}
j_l(x) = (-x)^l \left( \frac{1}{x} \frac{d}{dx} \right)^l \frac{\sin x}{x}
\quad \text{and} \quad
y_l(x) = -(-x)^l \left( \frac{1}{x} \frac{d}{dx} \right)^l \frac{\cos x}{x}.
\end{align}
For integer $l$, the spherical Bessel functions of the first and second kind are related to each other by
\begin{align}
j_{-(l+1)}(x) = \sqrt{\frac{\pi}{2x}} J_{-l - 1/2}(x) = (-1)^{l + 1} y_l(x).
\end{align}

In this paper, we concern ourselves only with the spherical Bessel functions of the first kind with integer order $l \ge 0$. The first few $j_l(x)$ are as follows.
\begin{subequations}
\begin{align}
j_0(x) &= \frac{\sin x}{x} = \mathrm{sinc}(x)
\\
j_1(x) &= \frac{\sin x}{x^2} - \frac{\cos x}{x}
\\
j_2(x) &= \left(\frac{3}{x^2} - 1\right) \frac{\sin x}{x} - \frac{3 \cos x}{x^2}
\end{align}
\end{subequations}
This pattern continues for higher orders, and all $j_l(x)$ with natural $l$ can be written as a sum of sines and cosines multiplied by inverse powers of $x$.

We will restrict ourselves to positive arguments $x$, although $j_l(x)$ is defined over the entire complex plane. Negative arguments of $j_l(x)$ are related to positive arguments by
\begin{align} \label{eq:parity}
j_l(-x) = (-1)^l j_l(x).
\end{align}
(A more general relation applies for complex $x$ or non-integer $l$.) 

Spherical Bessel functions satisfy two important recursion relations which can be combined to construct further relations.
\begin{subequations} \label{eq:recursions}
\begin{align}
j_l(x) &= \frac{l-1}{x} j_{l-1}(x) - \partial_x j_{l-1}(x)
\\
j_l(x) &= \frac{2l - 1}{x} j_{l-1}(x) - j_{l-2}(x)
\end{align}
\end{subequations}

Spherical Bessel functions satisfy a closure relation
\begin{align}
\int_0^\infty x^2 j_l(kx) j_l(k' x) dx = \frac{\pi}{2 k^2} \delta(k - k')
\end{align}
where $\delta$ is the Dirac delta function, as well as an orthogonality relation
\begin{align}
\int_{-\infty}^\infty j_k(x) j_l(x) dx = \frac{\pi}{2l+1} \delta_{kl}
\end{align}
for $k, l \in \mathbb{N}$, where $\delta_{kl}$ is the Kronecker delta.

A number of infinite integrals over spherical Bessel functions are known \cite{Mehrem2009}.
\begin{align}
\int_0^\infty j_l(x) dx &= \frac{\sqrt{\pi} \, \Gamma \! \left(\frac{l+1}{2}\right)}{2 \, \Gamma \! \left(1 + \frac{l}{2}\right)}
\\
\int_0^\infty j_l(x)^2 dx &= \frac{\pi}{2(2l+1)}
\\
\int_0^\infty j_l(kx) j_l(k'x) dx &= \frac{\pi}{2(2l+1)} \frac{(k_<)^l}{(k_>)^{l+1}}
\end{align}
Here, $k_<$ and $k_>$ are the smaller and larger of $k$ and $k'$ respectively.

The recursion relationships described in the rest of this paper excel when performing integrals over many oscillations. If one is trying to perform integrals over few if any oscillations at low $x$, quadrature methods will likely be better than the recursion methods, which can develop formally canceling divergences at $x = 0$ (this is a well-known issue in evaluating Bessel function recursion relations \cite{Rosenheinrich2016}). From linear regression on $0 \le l \le 100$, the first zero of $j_l(x)$ occurs roughly around
\begin{align}
x \approx 4.75 + 1.05 l.
\end{align}
Knowledge of the location of the first zero of the integrand can help in determining which computational method to use. Note that this is an overestimate for very low $l$; $j_0(x)$ has its first zero at $x = \pi$. 

It is also worth knowing the behavior of $j_l(x)$ near $x = 0$. At $x=0$, $j_0(0) = 1$, while $j_l(x) = 0$ for $l > 0$. The series expansion of $j_l(x)$ about $x = 0$ is given by
\begin{align}
j_l(x) = x^l \left(\frac{\sqrt{\pi}}{2^{l+1} \, \Gamma \! \left(l + \frac{3}{2}\right)} + O(x^2)\right).
\end{align}

\section{Preliminaries} \label{sec:prelim}

Because spherical Bessel functions can be written in terms of sines and cosines divided by powers of their argument, we will often encounter integrals of the form
\begin{align}
X_n(x) &= \int dx \, x^n \sin(x)
\\
Y_n(x) &= \int dx \, x^n \cos(x)
\end{align}
for integer $n$. It will be useful to be able to compute these efficiently. 

Integrating by parts, we obtain the following two recursion relations.
\begin{align}
X_n(x) &= n Y_{n-1}(x) - x^n \cos(x)
\\
Y_n(x) &= x^n \sin(x) - n X_{n-1}(x)
\end{align}
If $n > 0$, these relations allow us to step down by a single $n$, terminating at $n=0$ with
\begin{align}
X_0(x) = - \cos(x) \quad \text{and} \quad Y_0(x) = \sin(x).
\end{align}
For $n < 0$, we can invert these recursion relations to allow us to step upwards in $n$.
\begin{align}
X_n(x) &= \frac{1}{n+1} [x^{n+1} \sin(x) - Y_{n+1}(x)]
\\
Y_n(x) &= \frac{1}{n+1} [x^{n+1} \cos(x) + X_{n+1}(x)]
\end{align}
These recursions terminate at $n=-1$ with the special functions
\begin{align}
X_{-1}(x) &= \si(x) = \int_0^x \frac{\sin(t)}{t} dt
\\
Y_{-1}(x) &= \ci(x) = - \int_x^{\infty} \frac{\cos(t)}{t} dt
\end{align}
where $\si(x)$ is the sine integral, and $\ci(x)$ is the cosine integral. These special functions are available in scientific libraries such as SciPy and GSL.

We now look closely at $x=0$. We have $X_0(0) = -1$ and $Y_0(0) = 0$. These can be used to obtain $X_n(0) = Y_n(0) = 0$ for all $n > 0$. Looking towards negative $n$, we have $X_{-1}(0) = 0$, and $X_n(0)$ is divergent for $n < -1$. $Y_n(0)$ is divergent for all $n < 0$.

In the following sections, we will encounter expressions of the form $X_n(\alpha x)/\alpha^{n+1}$ and $Y_n(\alpha x)/\alpha^{n+1}$. In the regime where $\alpha x \lesssim \pi$, we recommend using numerical quadrature. For $n \ge 0$, one can also turn to series expansions.
\begin{align}
\frac{X_n(\alpha x)}{\alpha^{n+1}} &= \int dx \, x^n \sin(\alpha x)
=
- \frac{\Gamma(n+1) \cos(n \pi/2)}{\alpha^{n+1}} + \alpha x^{n+2} \left(\frac{1}{2+n} - \frac{\alpha^2 x^2}{6(4+n)} + \frac{\alpha^4 x^4}{120(6+n)} + O((\alpha x)^6)\right)
\\
\frac{Y_n(\alpha x)}{\alpha^{n+1}} &= \int dx \, x^n \cos(\alpha x)
=
\frac{\Gamma(n+1) \sin(n \pi/2)}{\alpha^{n+1}} + x^{n+1} \left(\frac{1}{1+n} - \frac{\alpha^2 x^2}{2(3+n)} + \frac{\alpha^4 x^4}{24(5+n)} + O((\alpha x)^6)\right)
\end{align}
Note the presence of constant terms that turn on and off with $n$. In the case of $n = -1$, we suggest using expansions for the sine and cosine integrals. For lower $n$, we were unable to find appropriate expansions, and urge care if numerical results are required in this regime.

\section{Integrals of a Single Spherical Bessel}

We now turn towards integrals of the form
\begin{align}
I^n_l(x; \alpha) = \int dx \, x^n j_l(\alpha x)
\end{align}
with $n \in \mathbb{Z}$, $l \in \mathbb{N}$, and $\alpha \in \mathbb{R}$. The special case of $\alpha = 0$ simplifies the integral drastically, and so we assume $\alpha \neq 0$. A substitution then yields
\begin{align}
I^n_l(x; \alpha) = \alpha^{-n-1} I^n_l(\alpha x)
\end{align}
where we define
\begin{align}
I^n_l(x) = \int dx \, x^n j_l(x).
\end{align}
Using the small $x$ behavior of $j_l(x)$, we see that $I^n_l(x)$ is finite at $x = 0$ for $l + n > -1$.

To evaluate $I^n_l(x)$, we begin by looking at $l = 0$. Because $j_0(x) = \sin(x) / x$, we immediately have
\begin{align}
I^n_0(x) = \int dx \, x^{n-1} \sin(x) = X_{n-1}(x).
\end{align}
Of particular note is $I^0_0(x) = \si(x)$.

For $l>0$, we use the recursion relationship
\begin{align}
j_l(x) &= \frac{l-1}{x} j_{l-1}(x) - \partial_x j_{l-1}(x)
\end{align}
to compute
\begin{align}
I^n_l(x) = \int dx \, x^n \left[\frac{l-1}{x} j_{l-1}(x) - \partial_x j_{l-1}(x)\right]
= (l + n - 1) I^{n-1}_{l-1}(x) - x^n j_{l-1}(x). \label{eq:Irecursion}
\end{align}
This allows $l$ to decrease by one each step, recursing until $l=0$ with $I^{n-l}_0(x) \equiv X_{n-l-1}(x)$.

The recursion will run through $(l,n)$ pairs as $(l-i,n-i)$ for integer $i$ where $0 \le i \le l$. There are $l+1$ steps, assuming the recursion doesn't truncate from the $(l + n - 1)$ term vanishing. To truncate, we must have $l + n$ odd. The recursion will then truncate when $i = (l + n - 1)/2$. Due to the range of $i$, this further requires $1 - l \le n \le 1 + l$.

Mathematica will integrate $I^n_l(x)$ directly, expressing the result in terms of a hypergeometric $\tensor[_1]{F}{_2}$ function. For special choices of parameters, the hypergeometric function yields simpler functions. These results can be directly verified by differentiation and use of the spherical Bessel recursion relations \eqref{eq:recursions}.
\begin{align}
I^{2+l}_l(x) &= x^{2+l} j_{l+1}(x) \label{eq:I1}
\\
I^{4+l}_l(x) &= x^{3+l} [(2l+3) j_{l+2}(x) - x j_{l+3}(x)] \label{eq:I2}
\\
I^{1-l}_l(x) &= \frac{\sqrt{\pi} 2^{-l}}{\Gamma \! \left(l+\frac{1}{2}\right)} - x^{1-l} j_{l-1}(x) \label{eq:I3}
\end{align}
Generally speaking, $n = l + \Omega$ and $n = 1 - l + \Omega$ for positive even $\Omega$ yield results similar to these, but the expressions quickly become more complicated as $\Omega$ increases.

\section{Integrals of a Spherical Bessel Squared}

Let us now look towards the following integral.
\begin{align}
H^n_l(x; \alpha) = \int dx \, x^n j_l(\alpha x)^2
\end{align}
Again, $\alpha = 0$ drastically simplifies the integral, and so we assume $\alpha \neq 0$. We again perform a substitution to obtain
\begin{align}
H^n_l(x; \alpha) = \alpha^{-n-1} H^n_l(\alpha x)
\end{align}
where we have defined
\begin{align}
H^n_l(x) = \int dx \, x^n j_l(x)^2.
\end{align}
From a series expansion of the integrand, $H^n_l(x)$ is finite at $x = 0$ for $2l + n > -1$.

We once again start with the $l=0$ case, writing $j_0(x) = \sin(x) / x$. By further writing $\sin^2(x) = (1-\cos(2x))/2$, we obtain
\begin{align}
H^n_0(x) = \int dx \, x^n j_0(x)^2
= \frac{1}{2} \int dx \, x^{n-2} - \frac{1}{2^n} Y_{n-2}(2x).
\end{align}
If $n \neq 1$, this becomes
\begin{align}
H^n_0(x) 
= \frac{x^{n-1}}{2(n-1)} - \frac{1}{2^n} Y_{n-2}(2x),
\end{align}
while for $n = 1$, we instead have
\begin{align}
H^1_0(x) = \frac{1}{2} \left[\ln(x) - \ci(2x)\right].
\end{align}
This is finite in the limit of small $x$,
\begin{align}
H^1_0(x) = -\frac{\gamma + \ln(2)}{2} + \frac{x^2}{2} - \frac{x^4}{12} + O(x^6)
\end{align}
where $\gamma$ is the Euler-Mascheroni constant.

We now turn to the case of $l > 0$. Note the following relationship, which can be verified by using the recursion relationships \eqref{eq:recursions}.
\begin{align}
\partial_x \left[\left(1 - \frac{n}{2}\right) x^{n-1} j_{l-1}(x)^2 - x^n j_{l-1}(x) j_l(x)\right]
=
\frac{1}{2} (2 - n) (2l + n - 3) x^{n-2} j_{l-1}(x)^2
+ x^n j_l(x)^2
- x^n j_{l-1}(x)^2
\end{align}
Performing an indefinite integral over $x$ leads to
\begin{align}
\left(1 - \frac{n}{2}\right) x^{n-1} j_{l-1}(x)^2 - x^n j_{l-1}(x) j_l(x)
=
\frac{1}{2} (2 - n) (2l + n - 3) H^{n-2}_{l-1}(x)
+ H^n_l(x)
- H^n_{l-1}(x)
\end{align}
where constants of integration are included in the $H$ functions. This result can be rearranged to solve for $H^n_l(x)$.
\begin{align} \label{eq:Hrecursion}
H^n_l(x) = H^n_{l-1}(x)
+ \frac{1}{2} (n - 2) (2l + n - 3) H^{n-2}_{l-1}(x)
+ \left(1 - \frac{n}{2}\right) x^{n-1} j_{l-1}(x)^2 - x^n j_{l-1}(x) j_l(x)
\end{align}
This allows us to recurse to lower $l$ values, which can be repeated until we get to $H^n_0(x)$, which has been evaluated above.

Mathematica will compute $H^n_l(x)$ directly in terms of a hypergeometric $\tensor[_2]{F}{_3}$ function. For special values of the parameters, the hypergeometric function simplifies to simpler functions. Once again, these results can be directly verified by differentiation and use of the spherical Bessel recursion relations \eqref{eq:recursions}.
\begin{align}
H^{(-1)}_l(x) &= \int dx \, \frac{j_l(x)^2}{x} = \frac{x^2 j_{l-1}(x)^2 - 2lx j_{l-1}(x) j_l(x) + (x^2-l) j_l(x)^2}{2l(l+1)} \label{eq:H1}
\\
H^{(1-2l)}_l(x) &= \int dx \, x^{1-2l} j_l(x)^2 = 
\frac{\pi}{4^{l+1}l \, \Gamma \! \left(l + \frac{1}{2}\right)^2}
-\frac{x^{2(1-l)} (j_{l-1}(x)^2 + j_l(x)^2)}{4 l} \label{eq:H2}
\\
H^2_l(x) &= \int dx \, x^2 j_l(x)^2 = \frac{x^3}{2} \left( j_l(x)^2 - j_{l-1}(x) j_{l+1}(x) \right) \label{eq:H3}
\\
H^4_l(x) &= \int dx \, x^4 j_l(x)^2
=
\frac{x^2}{12}\bigg(
-(2l+3)(4l^2+2x^2-1) j_l(x) j_{l-1}(x)
+x (4l(l+1)+2x^2-3) j_{l-1}(x)^2
\nonumber\\
&\qquad \qquad \qquad \qquad \qquad \qquad
+x (4l(l+2) + 2x^2 + 3) j_l(x)^2
\bigg) \label{eq:H4}
\\
H^{2l+3}_l(x) &= \int dx \, x^{2l+3} j_l(x)^2 = \frac{x^{2(l+2)}}{4(l+1)} \left(j_l(x)^2 + j_{l+1}(x)^2\right) \label{eq:H5}
\end{align}
Further closed-form expressions exist for $n = \Omega$, $n = 1 - 2l + \Omega$ and $n = 2l + 1 + \Omega$ with positive even $\Omega$, although the complexity rapidly increases with $\Omega$. It is useful to note that these analytic results often allow for an early termination to the recursion relationships, without having to recurse all the way down to $l = 0$.

\section{Integrals of Two Spherical Bessels, Same Order}

We now turn to integrals of two spherical functions that have the same order.
\begin{align}
K^n_l(x; \alpha, \beta) = \int dx \, x^n j_l(\alpha x) j_l(\beta x)
\end{align}
If $\alpha$ or $\beta = 0$, the integral simplifies significantly, while if $\alpha = \beta$, we reduce to a previous case. Hence, we assume $\alpha \neq \beta \neq 0$. Note that $K^n_l$ is symmetric in $\alpha$ and $\beta$. As for $H^n_l(x)$, $K^n_l(x; \alpha, \beta)$ is finite at $x = 0$ for $2l + n > -1$.

We once again begin with $l=0$. Writing $j_0(x) = \sin(x) / x$, we obtain
\begin{align}
K^n_0(x; \alpha, \beta) = \frac{1}{\alpha \beta} \int dx \, x^{n-2} \sin(\alpha x) \sin(\beta x).
\end{align}
Noting that 
\begin{align}
\sin(\alpha x) \sin(\beta x) = \frac{1}{2} \left[ \cos((\alpha - \beta) x) - \cos((\alpha + \beta) x)\right],
\end{align}
the integral can be evaluated to yield
\begin{align}
K^n_0(x; \alpha, \beta) = \frac{1}{2 \alpha \beta} \left[\frac{1}{(\alpha - \beta)^{n-1}} Y_{n-2}((\alpha - \beta) x) - \frac{1}{(\alpha + \beta)^{n-1}} Y_{n-2}((\alpha + \beta) x)\right].
\end{align}
If $|\alpha| \sim |\beta|$, then care must be taken in computing these quantities (see Section \ref{sec:prelim}).

For $l > 0$, we again construct a derivative relationship.
\begin{align}
&\partial_x \left[\frac{x^{n-1}}{2 \alpha \beta}\left[(2 - n) j_{l - 1}(\alpha x) j_{l - 1}(\beta x) - x (\beta j_{l - 1}(\alpha x) j_l(\beta x) + \alpha j_l(\alpha x) j_{l - 1}(\beta x))\right] \right]
\nonumber\\
&\quad = x^n j_l(\alpha x) j_l(\beta x) - \frac{1}{2 \alpha \beta} (n - 2) (n + 2l - 3) x^{n-2} j_{l-1}(\alpha x) j_{l-1}(\beta x)
- \frac{\alpha^2 + \beta^2}{2 \alpha \beta} x^n j_{l-1}(\alpha x) j_{l-1}(\beta x)
\end{align}
This can be derived using the recursion relations \eqref{eq:recursions}. Performing an indefinite integral over $x$, we obtain
\begin{align}
&\frac{x^{n-1}}{2 \alpha \beta}\left[(2 - n) j_{l - 1}(\alpha x) j_{l - 1}(\beta x) - x (\beta j_{l - 1}(\alpha x) j_l(\beta x) + \alpha j_l(\alpha x) j_{l - 1}(\beta x))\right]
\nonumber\\
&\quad = K^n_l(x;\alpha,\beta) - \frac{1}{2 \alpha \beta} (n - 2) (n + 2l - 3) K^{n-2}_{l-1}(x;\alpha,\beta)
- \frac{\alpha^2 + \beta^2}{2 \alpha \beta} K^n_{l-1}(x;\alpha,\beta) \,.
\end{align}
Once again, constants of integration are absorbed into the definition of $K$. Solving for $K^n_l(x;\alpha, \beta)$, we obtain the following recursion relation.
\begin{align} \label{eq:Krecursion}
K^n_l(x;\alpha,\beta) &=
\frac{1}{2 \alpha \beta} \bigg[
(\alpha^2 + \beta^2) K^n_{l-1}(x;\alpha,\beta)
+ (n - 2) (n + 2l - 3) K^{n-2}_{l-1}(x;\alpha,\beta)
\nonumber\\
&\qquad+
(2 - n) x^{n-1} j_{l - 1}(\alpha x) j_{l - 1}(\beta x) 
- x^n [\beta j_{l - 1}(\alpha x) j_l(\beta x) + \alpha j_l(\alpha x) j_{l - 1}(\beta x)] \bigg]
\end{align}
In the limit $\alpha = \beta = 1$, this reduces to the recursion relation \eqref{eq:Hrecursion} for $H^n_l(x)$.

The only closed-form solution we could find for integrals of this type is
\begin{align}
K^2_l(x; \alpha, \beta) = \int dx \, x^2 j_l(\alpha x) j_l(\beta x) = 
\frac{x^2}{\alpha ^2 - \beta^2} \left[\beta j_l(\alpha x) j_{l-1}(\beta x) - \alpha j_{l-1}(\alpha x) j_l(\beta x)\right]. \label{eq:K1}
\end{align}
Once again, this can be checked by differentiation and application of the recursion relations \eqref{eq:recursions}.

\section{Integrals of Two Spherical Bessels, Different Order}

We now consider our final integral, where we allow for the order of the Bessel functions to be different.
\begin{align}
L^n_{kl}(x; \alpha, \beta) = \int dx \, x^n j_k(\alpha x) j_l(\beta x)
\end{align}
If $\alpha = 0$ or $\beta = 0$, the integral once again simplifies drastically, and so we assume that $\alpha$ and $\beta$ are nonzero real numbers. From series expansions of the integrand, $L^n_{kl}(x; \alpha, \beta)$ is finite at $x = 0$ for $k + l + n > -1$.

Our approach in this case is to reduce the larger of the two orders through repeated application of the recursion relations \eqref{eq:recursions}. Assuming that $k < l$, we obtain
\begin{align} \label{eq:Lrecurse}
L^n_{kl}(x; \alpha, \beta) = \int dx \, x^n j_k(\alpha x) \left[\frac{2l - 1}{\beta x} j_{l-1}(\beta x) - j_{l-2}(\beta x)\right]
= \frac{2l - 1}{\beta} L^{n-1}_{k(l-1)}(x; \alpha, \beta) - L^n_{k(l-2)}(x; \alpha, \beta).
\end{align}
We can use this repeatedly until either $l = k$, in which case the $K$ integrals apply, or $l = k + 1$, which we now evaluate explicitly.

Consider the result
\begin{align}
\partial_x \big(x^{n+1} j_{l-1}(\alpha x) j_l(\beta x)\big) = 
(n-1) x^n j_{l-1}(\alpha x) j_l(\beta x)
+ \beta x^{n+1} j_{l-1}(\alpha x) j_{l-1}(\beta x) - \alpha x^{n+1} j_l(\alpha x) j_l(\beta x)
\end{align}
which can be verified by differentiation and application of the recursion relations \eqref{eq:recursions}. Integrating this over $x$ yields
\begin{align}
x^{n+1} j_{l-1}(\alpha x) j_l(\beta x) = 
(n-1) L^n_{l-1,l}(x; \alpha, \beta)
+ \beta K^{n+1}_{l-1}(x; \alpha, \beta)
- \alpha K^{n+1}_l(x; \alpha, \beta)
\end{align}
where constants of integration are absorbed into the functions. Rearranging, we then find
\begin{align} \label{eq:Lclose}
L^n_{l-1,l}(x; \alpha, \beta)
=
\frac{1}{n-1}\big(
x^{n+1} j_{l-1}(\alpha x) j_l(\beta x) 
+ \alpha K^{n+1}_l(x; \alpha, \beta)
- \beta K^{n+1}_{l-1}(x; \alpha, \beta)
\big).
\end{align}

The recursion relations may also terminate at $L^n_{10}(x; \alpha, \beta) = L^n_{01}(x; \beta, \alpha)$, which we now compute.
\begin{align}
L^n_{01}(x; \alpha, \beta) &= \int dx \, x^n \frac{\sin(\alpha x)}{\alpha x} \left(\frac{\sin(\beta x)}{\beta^2 x^2} - \frac{\cos(\beta x)}{\beta x}\right)
\\
&= \frac{1}{2 \alpha \beta^2}\int dx \, x^{n-3} \bigg[\cos((\alpha - \beta)x) - \cos((\alpha + \beta)x)\bigg]
\nonumber\\
& \qquad \qquad \qquad - \frac{1}{2 \alpha \beta}\int dx \, x^{n-2} \bigg[\sin((\alpha - \beta)x) + \sin((\alpha + \beta)x)\bigg]
\\
&= 
\frac{1}{2 \alpha \beta^2} \bigg[\frac{1}{(\alpha - \beta)^{n-2}} Y_{n-3}((\alpha - \beta)x) - \frac{1}{(\alpha + \beta)^{n-2}} Y_{n-3}((\alpha + \beta)x)\bigg]
\nonumber\\
& \qquad \qquad \qquad - \frac{1}{2 \alpha \beta} \bigg[\frac{1}{(\alpha - \beta)^{n-1}} X_{n-2}((\alpha - \beta)x) + \frac{1}{(\alpha + \beta)^{n-1}} X_{n-2}((\alpha + \beta)x)\bigg]
\end{align}
If $|\alpha| \sim |\beta|$, then these quantities will again need to be computed carefully (see Section \ref{sec:prelim}).

We were only able to find one closed-form relation for general $L^n_{kl}(x;\alpha,\beta)$. 
\begin{align}
&(\alpha^2 - \beta^2) L^2_{kl}(x; \alpha, \beta) + [l(l+1) - k(k+1)] L^0_{kl}(x; \alpha, \beta) 
\nonumber\\
&\qquad=
\int dx \, [(\alpha^2 - \beta^2) x^2 + l(l+1) - k(k+1)] j_k(\alpha x) j_l (\beta x)
\\
&\qquad=
\beta x^2 j_k(\alpha x) j_{l-1}(\beta x) - \alpha x^2 j_{k-1}(\alpha x) j_l(\beta x) 
+ (k-l) x j_k(\alpha x) j_l(\beta x)
\label{eq:Lresult}
\end{align}
Unfortunately, this is not a particularly useful formula unless the exact combination depicted here arises. However, special limits such as $\alpha = \beta$ or $k = l$ do yield useful results, as shown above and below.

If we restrict ourselves to the case of $\alpha = \beta$, much more can be derived. Let us begin by defining
\begin{align}
L^n_{kl}(x) = \int dx \, x^n j_k(x) j_l(x)
\end{align}
so that
\begin{align}
L^n_{kl}(x;\alpha,\alpha) = \alpha^{-1-n} L^n_{kl}(\alpha x).
\end{align}
Note that $L^n_{kl}(x)$ is symmetric under interchange of $k$ and $l$. Again, the basic idea will be to use a recursion relation to decrease the highest of $k$ or $l$ by using Eq. \eqref{eq:Lrecurse}. In the case of $k = l - 1$, instead of Eq. \eqref{eq:Lclose}, we have the simpler result
\begin{align} \label{eq:Lrecurse2}
L^n_{(l-1)l}(x) = (l-1) H^{n-1}_{l-1}(x) - \int dx \, x^n j_{l-1}(x) \partial_x j_{l-1}(x)
= \left(l+\frac{n}{2}-1\right) H^{n-1}_{l-1}(x)
- \frac{1}{2} x^n j_{l-1}(x)^2.
\end{align}
If needed, $L^n_{01}(x)$ can be computed by
\begin{align}
L^n_{01}(x) &= \int dx \, x^{n-2} \sin(x) \left(\frac{\sin (x)}{x} - \cos (x)\right) 
\\
&= \int dx \, \left(\frac{x^{n-3}}{2} - \frac{x^{n-3} \cos(2x)}{2} - \frac{x^{n-2} \sin(2x)}{2}\right)
\\
&= \int dx \, \frac{x^{n-3}}{2} - \frac{1}{2^{n-1}} Y_{n-3}(2x) - \frac{1}{2^n} X_{n-2}(2x).
\end{align}

$L^n_{kl}(x)$ can be expressed directly in terms of a rather messy hypergeometric $\tensor[_3]{F}{_4}$ function. We found that certain special choices of parameters lead to simplifications in the hypergeometric function that allow the result to be written quite cleanly. As always, these results can be verified by differentiation and application of recursion relations \eqref{eq:recursions}.
\begin{align}
L^0_{kl}(x) &= \int dx \, j_k(x) j_l(x)
= \frac{x}{k(k+1) - l(l+1)} \left[ 
x (j_{k-1}(x) j_l(x)
- j_k(x) j_{l-1}(x))
+ (l-k) j_k(x) j_l(x)
\right] \label{eq:L1}
\\
L^{-1}_{kl}(x) &= \int dx \, \frac{j_k(x) j_l(x)}{x}
= \frac{j_k(x) j_l(x)}{k + l}
- \frac{x j_{k+1}(x) j_l(x)}{(1 + k - l)(k+l)}
- \frac{x j_k(x) j_{l+1}(x)}{(1 + l - k)(k+l)}
\nonumber\\
&\qquad \qquad \qquad \qquad \qquad \qquad
+ \frac{2x^2 (j_{k+1}(x) j_{l+1}(x) + j_k(x) j_l(x))}{(k+l)(2+k+l)(1+k-l)(1+l-k)} \label{eq:L2}
\\
L^{1-k-l}_{kl}(x) &= \int dx \, x^{1-k-l} j_k(x) j_l(x)
= 
\frac{\pi}{2^{k + l + 1}(k+l) \, \Gamma \! \left(k+\frac{1}{2}\right) \Gamma \! \left(l+\frac{1}{2}\right)}
\nonumber\\
&\qquad \qquad \qquad \qquad \qquad \qquad \qquad 
-\frac{x^{2-k-l} (j_{k-1}(x) j_{l-1}(x) + j_k(x) j_l(x))}{2 (k+l)} \label{eq:L3}
\\
L^{l-k+2}_{kl}(x) &= \int dx \, x^{l-k+2} j_k(x) j_l(x)
= \frac{x^{-k+l+3}}{2 (k-l-1)} \bigg(j_{k-1}(x) j_{l+1}(x) - j_k(x) j_l(x)\bigg) \label{eq:L4}
\\
L^{k+l+3}_{kl}(x) &= \int dx \, x^{k+l+3} j_k(x) j_l(x)
= \frac{x^{k+l+4}}{2 (k+l+2)} \bigg(j_{k+1}(x) j_{l+1}(x) + j_k(x) j_l(x)\bigg) \label{eq:L5}
\end{align}
Further closed-form expressions are available for $n = k - l + \Omega$, $n = l - k + \Omega$, $n = k + l + 1 + \Omega$ and $n = 1 - k - l + \Omega$ with positive even $\Omega$, although once again, the expressions quickly become complicated as $\Omega$ increases.

\section{Conclusions}

We have investigated indefinite integrals of powers of $x$ multiplying one or two spherical Bessel functions of the first kind with integer order $l \ge 0$. In the process, we presented a number of closed-form indefinite integrals and constructed recursion relations to allow the numerical computation of all the integrals we considered. 

The recursion techniques that we have presented are highly complementary to numerical quadrature methods. When integrating spherical Bessel functions over many oscillations, our expressions allow for accurate computations in situations where quadrature methods suffer. On the other hand, our recursion relations have numerical issues at small $x$, where formal divergences are canceling, which ruins numerical precision. In this regime, there are few if any oscillations and quadrature methods excel.

In the regime in which we don't need to worry about canceling divergences, we have compared a numerical implementation of our results to Mathematica's high-precision numerical integration, with excellent results.

\acknowledgments

This work is supported in part by the U.S. Department of Energy under grant Contract Number DE-SC0012567, and in part by MIT's Undergraduate Research Opportunities Program.

\appendix

\section{Results translated into ordinary Bessel functions} \label{appendix}

In this appendix, we translate our results into expressions involving ordinary Bessel functions. While the results translate perfectly well, the recursion relations typically do not terminate. Apart from the recursion relations, a number of these results appear in the literature, e.g., \cite{wolfram-bessel}.

\subsection{Single Bessel integrals}

The following are translations of Eqs. \eqref{eq:Irecursion}, \eqref{eq:I1}, \eqref{eq:I2} and \eqref{eq:I3} respectively.
\begin{align}
\int dx \, x^n J_{l}(x) &= (l + n - 1) \int dx \, x^{n-1} J_{l - 1}(x) - x^n J_{l - 1}(x)
\\
\int dx \, x^{l+1} J_{l}(x) &= x^{l+1} J_{l+1}(x)
\\
\int dx \, x^{l+3} J_{l}(x) &= x^{2+l} [2(l+1) J_{l + 2}(x) - x J_{l + 3}(x)]
\\
\int dx \, x^{1 - l} J_{l}(x) &= \frac{2^{-l}}{\Gamma \! \left(l\right)} - x^{1-l} J_{l - 1}(x)
\end{align}

\subsection{Bessel Squared integrals}

The following are translations of Eqs. \eqref{eq:Hrecursion}, \eqref{eq:H1}, \eqref{eq:H2}, \eqref{eq:H3}, \eqref{eq:H4} and \eqref{eq:H5} respectively.

\begin{align}
\int dx \, x^n J_{}(x)^2
&= \int dx \, x^n J_{l - 1}(x)^2
+ \frac{1}{2} (n - 1) (2l + n - 3) \int dx \, x^{n-2} J_{l - 1}(x)^2
\nonumber
\\
&\qquad \qquad
+ \frac{1-n}{2} x^{n-1} J_{l - 1}(x)^2 - x^n J_{l}(x) J_{l - 1}(x)
\\
\int dx \, x^{-2} J_{l}(x)^2 &= \frac{2 x^2 J_{l - 1}(x)^2 - (4l-2)x J_{l - 1}(x) J_{l}(x) + (2 x^2 + 1 - 2l) J_{l}(x)^2}{x(4l^2-1)}
\\
\int dx \, x^{1-2l} J_{l}(x)^2 &= 
\frac{2}{4^l(2l-1) \Gamma(l)^2}
-\frac{x^{2(1-l)} (J_{l - 1}(x)^2 + J_{l}(x)^2)}{2 (2l-1)}
\\
\int dx \, x J_{l}(x)^2 &= \frac{x^2}{2} \left( J_{l}(x)^2 - J_{l - 1}(x) J_{l + 1}(x) \right)
\\
\int dx \, x^3 J_{l}(x)^2
&=
\frac{x}{6}\bigg(
-(l+1)(4l^2-4l+2x^2) J_{l}(x) J_{l - 1}(x)
+x (2l^2+x^2-2) J_{l - 1}(x)^2
\nonumber\\
&\qquad \qquad \qquad \qquad \qquad \qquad
+x (2l^2 + 2l + x^2) J_{l}(x)^2
\bigg)
\\
\int dx \, x^{2l+1} J_{l}(x)^2 &= \frac{x^{2l+2}}{2(2l+1)} \left(J_{l}(x)^2 + J_{l + 1}(x)^2\right)
\end{align}

\subsection{Two Bessels, Same Order}

The following are translations of Eqs. \eqref{eq:Krecursion} and \eqref{eq:K1} respectively.
\begin{align}
2 \alpha \beta \int dx \, x^{n} J_{l}(\alpha x) J_{l}(\beta x)
&=
(\alpha^2 + \beta^2) \int dx \, x^{n} J_{l - 1}(\alpha x) J_{l - 1}(\beta x)
\nonumber\\
&\qquad+
(n - 1) (n + 2l - 3) \int dx \, x^{n-2} J_{l-1}(\alpha x) J_{l-1}(\beta x)
\nonumber\\
&\qquad+
(1 - n) x^{n-1} J_{l-1}(\alpha x) J_{l-1}(\beta x)
- x^{n} [\beta J_{l-1}(\alpha x) J_{l}(\beta x) + \alpha J_{l}(\alpha x) J_{l - 1}(\beta x)]
\\
\int dx \, x J_{l}(\alpha x) J_{l}(\beta x) &= 
\frac{x}{\alpha ^2 - \beta^2} \left[\beta J_{l}(\alpha x) J_{l - 1}(\beta x) - \alpha J_{l - 1}(\alpha x) J_{l}(\beta x)\right]
\end{align}

\subsection{Two Bessels, Different Order}

The following are translations of Eqs. \eqref{eq:Lrecurse}, \eqref{eq:Lclose} and \eqref{eq:Lresult} respectively.
\begin{align}
\int dx \, x^n J_k(\alpha x) J_l(\beta x)
&= \frac{2(l - 1)}{\beta} \int dx \, x^{n-1} J_k(\alpha x) J_{l-1}(\beta x)
-
\int dx \, x^n J_k(\alpha x) J_{l-2}(\beta x)
\\
n \int dx \, x^n J_{l-1}(\alpha x) J_l(\beta x)
&=
x^{n+1} J_{l-1}(\alpha x) J_l(\beta x) 
+ \alpha \int dx \, x^{n+1} J_l(\alpha x) J_l(\beta x)
\nonumber\\
& \qquad 
- \beta \int dx \, x^{n+1} J_{l-1}(\alpha x) J_{l-1}(\beta x)
\\
(\alpha^2 - \beta^2) \int dx \, x J_k(\alpha x) J_l (\beta x)
&=
(k^2 - l^2) \int dx \, \frac{J_k(\alpha x) J_l (\beta x)}{x}
+ \beta x J_k(\alpha x) J_{l-1}(\beta x) 
\nonumber \\
& \qquad \qquad \qquad 
- \alpha x J_{k-1}(\alpha x) J_l(\beta x) 
+ (k-l) J_k(\alpha x) J_l(\beta x)
\end{align}
The following are translations of Eqs. \eqref{eq:Lrecurse2}, \eqref{eq:L1}, \eqref{eq:L2}, \eqref{eq:L3}, \eqref{eq:L4} and \eqref{eq:L5} respectively.
\begin{align}
\int dx \, x^n J_{l-1}(x) J_l(x) 
&= \left(l+\frac{n}{2} - 1\right) \int dx \, x^{n-1} J_{l-1}(x)^2
- \frac{1}{2} x^n J_{l-1}(x)^2
\\
\int dx \, \frac{J_k(x) J_l(x)}{x}
&= \frac{1}{k^2 - l^2} \bigg(
x [J_{k-1}(x) J_l(x)
- J_k(x) J_{l-1}(x)]
+ (l-k) J_k(x) J_l(x)
\bigg)
\\
\int dx \, \frac{J_k(x) J_l(x)}{x^2}
&= \frac{J_k(x) J_l(x)}{x(k + l - 1)}
- \frac{J_{k+1}(x) J_l(x)}{(1 + k - l)(k+l - 1)}
- \frac{J_k(x) J_{l+1}(x)}{(1 + l - k)(k+l - 1)}
\nonumber\\
&\qquad \qquad \qquad \qquad
+ \frac{2x [J_{k+1}(x) J_{l+1}(x) + J_k(x) J_l(x)]}{(k+l - 1)(k+l + 1)(1+k-l)(1+l-k)}
\\
\int dx \, x^{1-k-l} J_k(x) J_l(x)
&= 
\frac{1}{2^{k + l - 1}(k+l-1) \Gamma(k) \Gamma(l)}
-\frac{x^{2-k-l} (J_{k-1}(x) J_{l-1}(x) + J_k(x) J_l(x))}{2 (k+l-1)}
\\
\int dx \, x^{l-k+1} J_{k}(x) J_{l}(x)
&= \frac{x^{l-k+2}}{2 (k-l-1)} \bigg(J_{k - 1}(x) J_{l + 1}(x) - J_{k}(x) J_{l}(x)\bigg)
\\
\int dx \, x^{k+l+1} J_{k}(x) J_{l}(x)
&= \frac{x^{k+l+2}}{2 (k+l+1)} \bigg(J_{k + 1}(x) J_{l + 1}(x) + J_{k}(x) J_{l}(x)\bigg)
\end{align}

\bibliographystyle{utphys}
\bibliography{citations}

\end{document}